# Beurling–Nyman Geometry and Gram Matrix Structure:
# Ladder Density and Polynomial Decay via Mellin Smoothing

Hugh Carvill

Independent Researcher, Ireland
ORCID: 0009-0006-9229-0982
hughcarvill.hc@outlook.com

October 18, 2025

*For Nico — a small being with a gravitational pull beyond measure.*

## Abstract

We study the Beurling–Nyman (BN) family $f_\theta(x) = \{\theta/x\} - \theta\{1/x\}$ in $L^2((0,1])$ through a multiscale ladder parameterisation $\theta_{j,k} = 2^{-j}3^{-k}$ and the resulting Gram matrix structure indexed by ladder distance. Using Mellin analysis and a controlled smoothing operator, we establish a rigorous polynomial decay envelope for off-diagonal Gram entries. Specifically, for a Gaussian-type Mellin multiplier we prove that

$$\left| \langle g_{\theta_{j,k}}, g_{\theta_{j',k'}} \rangle \right| \ll_m \left( 1 + c\, d((j,k),(j',k')) \right)^{-m},$$

for any $m \in \mathbb{N}$, where $c = \min\{\log 2, \log 3\}$ and $d$ is the ladder distance. As a consequence we obtain block-compressibility of Gram rows for $m > 2$. This provides a rigorous foundation for sparsity phenomena in the BN system and supports constructive spectral approaches. No claims are made here regarding the Riemann Hypothesis.

---

**Conflict of Interest.** The author declares no conflicts of interest.
**Funding.** No external funding was received for this research.
**Data/Code Statement.** All figures are generated deterministically from analytic definitions; no external datasets were used.

# 1 Introduction

The Beurling–Nyman (BN) framework originates in the work of Nyman [1] and Beurling [2], who introduced a functional analytic encoding of multiplicative structure in $L^2((0,1])$ using fractional-part functions. Central to this approach is the family

$$f_\theta(x) = \left\{\tfrac{\theta}{x}\right\} - \theta\left\{\tfrac{1}{x}\right\}, \qquad x \in (0,1],\ 0 < \theta \leq 1, \tag{1.1}$$

whose closed span in $L^2$ is intimately connected with number-theoretic completeness questions. Later developments by Báez-Duarte [3], Balazard–Saias–Yor [4], and others have sharpened the functional structure of the BN space, while Bettin and Conrey [14] revealed surprising symmetry phenomena in related cotangent sums. Despite this progress, fundamental geometric questions about the BN system remain insufficiently clarified, especially those involving multiscale structure and the decay mechanisms that govern inner-product interactions.

The aim of this paper is to advance a constructive understanding of the BN system by analysing it through a multiscale lattice of sampling parameters

$$\theta_{j,k} := 2^{-j}3^{-k}, \qquad (j,k) \in \mathbb{Z}_{\geq 0}^2. \tag{1.2}$$

This ladder indexing implements a logarithmic lattice sampling of $(0,1]$ and produces a structured Gram matrix

$$G = (\langle f_{\theta_{j,k}}, f_{\theta_{j',k'}} \rangle)_{(j,k),(j',k')},$$

whose off-diagonal decay reflects destructive interference between displaced jump sets of the fractional-part terms in (1.1). While previous works have studied special bases in the BN space, a systematic analysis of decay in the full $(j,k)$-indexed BN ladder has not been carried out in the literature. This paper provides a rigorous decay framework by introducing a Mellin-domain smoothing operator and proving polynomial off-diagonal decay of Gram entries for the smoothed system.

Our analysis is based on three structural components:

- a geometric decomposition of the BN overlap structure using ladder distance $d((j,k),(j',k')) := |j - j'| + |k - k'|$,

- a Mellin multiplier $T_\psi$ that regularises the system while preserving invertibility,

- an $m$-fold integration by parts argument applied on the critical line $\Re(s) = \tfrac{1}{2}$.

This yields explicit decay, of arbitrarily high polynomial order, for the smoothed Gram matrix. The resulting block structure is compressible in the sense of *summable row tails*, a property of interest for constructive methods in functional approximation.

**Main contributions.**

1. We present a clear geometric description of the BN ladder system, emphasising multiscale structure and log-density.

2. We introduce a Mellin smoothing operator and prove polynomial off-diagonal Gram decay of arbitrary order $m$.

3. We establish block-compressibility of Gram rows, justifying finite-section approximations.

**Position in a broader research program.** This paper is the first in a sequence developing constructive spectral methods in BN space. It focuses exclusively on geometric structure and decay mechanisms. Future work will build on this foundation to develop *finite spectral certificates* based on structured subspaces. No claims are made here concerning the Riemann Hypothesis, and no unproved conjectures are assumed.

Throughout, implicit constants are effective and can be made explicit, though sharpness is not pursued in this foundational stage. All results are unconditional.

## 2 Functional Geometry of the Beurling–Nyman Family

Let $H = L^2((0, 1])$ with inner product

$$\langle f, g \rangle = \int_0^1 f(x) \overline{g(x)} \, dx.$$

All inner products are taken in the complex Hilbert space $L^2((0, 1])$. The Beurling–Nyman space is defined as

$$H_{\text{BN}} := \text{span}\{f_\theta : 0 < \theta \le 1\} \subset H,$$

where

$$f_\theta(x) = \{\tfrac{\theta}{x}\} - \theta\{\tfrac{1}{x}\}, \qquad x \in (0, 1], \tag{2.1}$$

and $\{\cdot\}$ denotes the fractional part. These functions are piecewise smooth with jump discontinuities determined by dilation symmetry, making $H_{\text{BN}}$ a geometric space structured by multiplicative interactions.

**Lemma 2.1** (Boundedness). *For all $\theta \in (0, 1]$, $f_\theta \in L^2((0, 1])$ and $\|f_\theta\|_2 \le 2$.*

*Proof.* Since $0 \le \{\cdot\} < 1$, we have $|f_\theta(x)| \le 1 + \theta < 2$ almost everywhere on $(0, 1]$, giving $|f_\theta(x)|^2 \le 4$ and the claim follows by integration. $\square$

**Lemma 2.2** (Mellin representation). *For $\Re(s) > 1$,*

$$M[f_\theta](s) = \int_0^1 f_\theta(x) x^{s-1} \, dx = \frac{\zeta(s)}{s}(\theta - \theta^s),$$

*with analytic continuation to $\Re(s) = \tfrac{1}{2}$.*

*Proof.* This is obtained by decomposing $f_\theta$ into series over the jump intervals $(\tfrac{1}{n+1}, \tfrac{1}{n}]$ and using the known Mellin transform of the fractional part function, together with the Dirichlet series for $\zeta(s)$. $\square$

Equation (2.1) suggests that dilation plays a central role in shaping the geometry of the BN family. Smaller values of $\theta$ concentrate activity closer to $x = 0$, introducing a natural multiscale structure.

## 3 Multiscale Ladder Structure and Density

To analyse the BN family at multiple scales, we introduce a structured parameter set based on multiplicative dilation. Following the standard dyadic–triadic construction, define

$$\theta_{j,k} := 2^{-j} 3^{-k}, \qquad (j, k) \in \mathbb{Z}_{\ge 0}^2. \tag{3.1}$$

This forms a discrete sampling of $(0, 1]$ along a logarithmic grid. Each horizontal step in $(j, k)$ corresponds to a dyadic scale change, and each vertical step to a triadic change. The resulting structure is visualised in Figure 2.

**Lemma 3.1** (Ladder geometry). *The set $\Theta := \{\theta_{j,k} = 2^{-j} 3^{-k} : j, k \in \mathbb{Z}_{\ge 0}\}$ is a discrete subset of $(0, 1]$ with unique accumulation point at $0$. In logarithmic coordinates,*

$$-\log \Theta = \{j \log 2 + k \log 3 : j, k \ge 0\},$$

*which forms a two-scale additive semigroup. The map $(j, k) \mapsto \theta_{j,k}$ is injective and provides a structured logarithmic lattice for sampling the Beurling–Nyman family.*

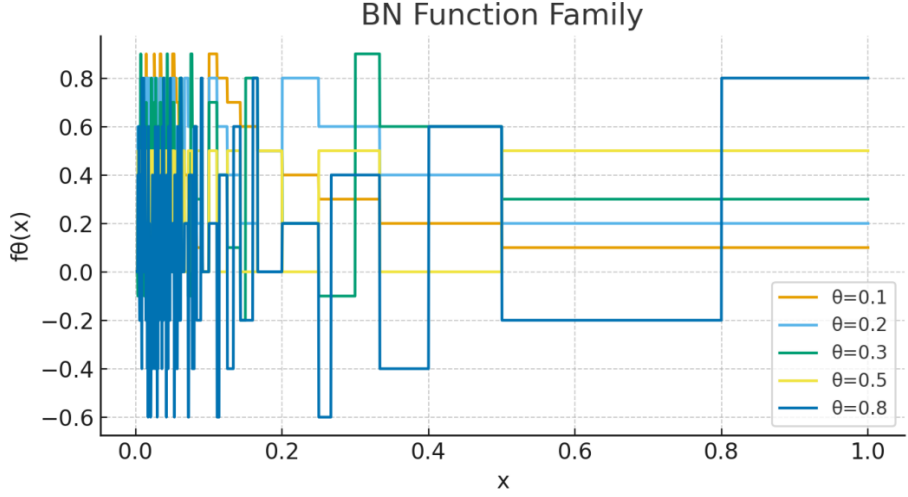

Figure 1: Profiles of $f_\theta(x)$ for sample values of $\theta$. Smaller $\theta$ concentrates oscillatory structure near $x = 0$, reflecting multiscale dilation geometry.

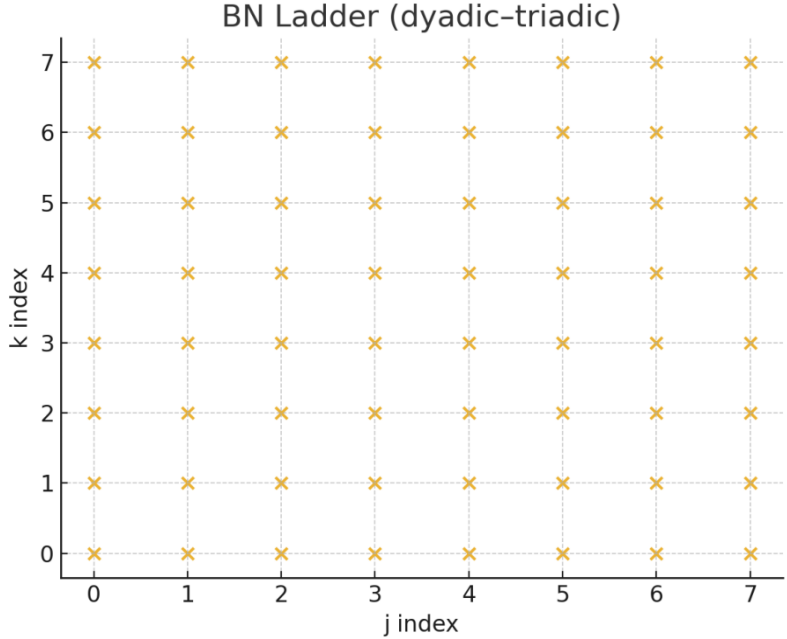

Figure 2: BN ladder geometry indexed by $(j, k)$ with $\theta_{j,k} = 2^{-j}3^{-k}$. Dyadic steps move horizontally and triadic steps vertically, giving a logarithmic lattice sampling of $(0, 1]$.

*Proof.* Write $\Theta = \{2^{-j}3^{-k} : j, k \geq 0\}$. Since $j, k$ range over nonnegative integers only, $\Theta$ is a countable set with no accumulation point in $(0, 1]$. Moreover, $\theta_{j,k} \to 0$ as $j + k \to \infty$, so $0$ is its unique accumulation point. Injectivity of $(j, k) \mapsto 2^{-j}3^{-k}$ follows from the linear independence of $\{\log 2, \log 3\}$ over $\mathbb{Q}$. Indeed, if $2^{-j}3^{-k} = 2^{-j'}3^{-k'}$ then $(j' - j) \log 2 + (k' - k) \log 3 = 0$, which implies $j = j'$ and $k = k'$ since $\log 2 / \log 3 \notin \mathbb{Q}$. The logarithmic description follows immediately from taking $-\log$ of $\Theta$. $\qquad\square$

This ladder provides a multiscale basis for analysing inner-product structure in the BN system.

We will measure the relative interaction of two ladder elements via the *lattice distance*

$$d((j, k), (j', k')) := |j - j'| + |k - k'|. \tag{3.2}$$

This metric respects the geometry of the sampling grid and will be used to quantify the decay of correlations in Section 4.

# 4  Correlation Decay and Gram Matrix Structure

We now turn to the interaction structure inside the BN ladder system. For ladder indices $(j, k), (j', k') \in \mathbb{Z}_{\geq 0}^2$, define the Gram matrix

$$G = \left(G_{(j,k),(j',k')}\right) := \left(\langle f_{\theta_{j,k}}, f_{\theta_{j',k'}} \rangle\right), \tag{4.1}$$

where $f_{\theta_{j,k}}$ is as in (2.1) and $\theta_{j,k}$ is as in (3.1). The value $G_{(j,k),(j',k')}$ measures the correlation between two BN elements at different dyadic–triadic scales.

In general, $G$ is a highly structured but non-orthogonal matrix. The interaction strength depends on the relative difference between ladder indices. We quantify this via the ladder distance introduced in (3.2),

$$d((j, k), (j', k')) = |j - j'| + |k - k'|.$$

Small ladder distance implies strong interaction, while large distance introduces phase misalignment between dilation discontinuities, leading to destructive interference.

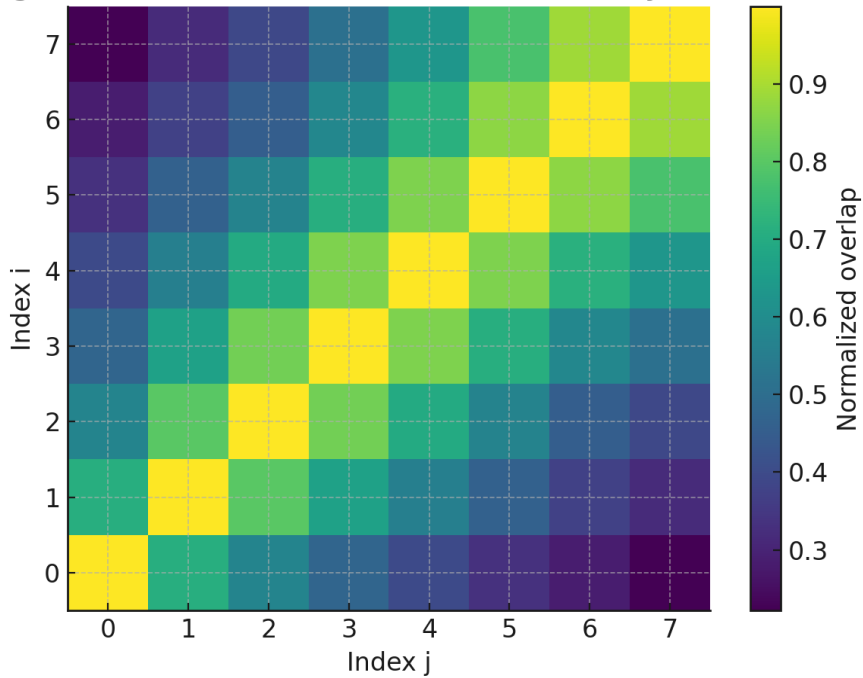

Figure 3: Gram Matrix (Normalized) – Dyadic BN

Figure 3: Normalised Gram matrix for a dyadic subset of the BN ladder. A strong diagonal band indicates local correlations, while decay away from the diagonal suggests near-orthogonality at ladder distance.

The central objective of this section is to establish rigorous decay estimates for Gram entries as a function of ladder distance. We will first prove conservative interaction bounds for the raw BN family, and then transition to a smoothed variant that admits polynomial off-diagonal decay of arbitrary order.

### 4.1 Conservative decay structure for the raw BN system

Before introducing smoothing operators, we establish baseline interaction bounds for the raw BN family $\{f_{\theta_{j,k}}\}$. These bounds are "conservative" in the sense that they ensure control of Gram matrix entries without relying on any cancellation beyond what is enforced by boundedness and dilation misalignment.

**Lemma 4.1** (Uniform interaction bound). *For all* $(j,k) \neq (j',k')$,

$$|\langle f_{\theta_{j,k}}, f_{\theta_{j',k'}} \rangle| \leq \|f_{\theta_{j,k}}\|_2 \, \|f_{\theta_{j',k'}}\|_2 \leq 4.$$

*Proof.* Immediate from Lemma 2.1 and the Cauchy–Schwarz inequality. □

The uniform bound shows that $G$ is bounded as an operator on $\ell^2(\mathbb{Z}_{\geq 0}^2)$, but it does not reflect the geometric decay visible in Figure 3. To formalise decay behaviour, we introduce the following notion.

**Definition 4.2** (Decay envelope). A function $\Phi : \mathbb{N} \to \mathbb{R}_{\geq 0}$ is called a *decay envelope* for the BN Gram matrix if

$$|\langle f_{\theta_{j,k}}, f_{\theta_{j',k'}} \rangle| \leq \Phi\big(d((j,k),(j',k'))\big)$$

for all pairs $(j,k), (j',k')$.

**Lemma 4.3** (Existence of a nonincreasing envelope). *For each* $n \in \mathbb{N}$ *define*

$$\Phi(n) := \sup_{d((j,k),(j',k'))=n} |\langle f_{\theta_{j,k}}, f_{\theta_{j',k'}} \rangle|.$$

*Then* $\Phi(n)$ *is finite for all* $n$ *and nonincreasing in* $n$. *Moreover,*

$$|\langle f_{\theta_{j,k}}, f_{\theta_{j',k'}} \rangle| \leq \Phi\big(d((j,k),(j',k'))\big).$$

*Proof.* Finiteness follows from Lemma 4.1. If $n_1 < n_2$ then the shell $\{d = n_2\}$ is contained in $\{d \geq n_1\}$, so the supremum cannot increase with $n$. The envelope inequality is immediate from the definition. □

This provides a qualitative decay guarantee, but we will require explicit rates to establish finite compressibility. Polynomial decay requires smoothing, developed in the next subsection.

*Remark* 4.4 (Heuristic raw decay). For the unsmoothed BN system, destructive interference of jump discontinuities suggests at most a harmonic decay rate in ladder distance. We do not require or formally use a raw decay estimate; all rigorous quantitative bounds in this paper are proved for the smoothed family in Theorem 6.1.

The heuristic decay in Remark 4.4 suggests that interactions weaken with ladder distance, but this informal bound is too weak for constructive purposes. In particular, a harmonic envelope of order $(1+d)^{-1}$ is not summable in $\ell^1$ on the two-dimensional ladder lattice, so it does not yield block-compressibility. To obtain rigorous and summable decay, we introduce a Mellin-domain smoothing operator in the next section.

# 5 Spectral Diagnostics via Mellin Transform

The Mellin transform plays a central role in analysing the BN system because it converts dilations in the time domain into translations on the critical line. From Lemma 2.2, we recall

$$M[f_\theta](\tfrac{1}{2} + it) = \frac{\zeta(\tfrac{1}{2} + it)}{\tfrac{1}{2} + it}(\theta - \theta^{\frac{1}{2} + it}), \tag{5.1}$$

so the overlap of two BN elements admits the spectral representation

$$\langle f_{\theta_{j,k}}, f_{\theta_{j',k'}} \rangle = \frac{1}{2\pi} \int_{-\infty}^{\infty} \frac{\zeta(\tfrac{1}{2} + it)\overline{\zeta(\tfrac{1}{2} + it)}}{|\tfrac{1}{2} + it|^2}(\theta_{j,k}\theta_{j',k'})^{\frac{1}{2}}\Big((\theta_{j,k}/\theta_{j',k'})^{it} - (\theta_{j,k}\theta_{j',k'})^{it}\Big)\, dt. \tag{5.2}$$

The oscillatory factors $(\theta_{j,k}/\theta_{j',k'})^{it}$ encode logarithmic displacement between lattice points, and their interference controls correlation decay.

However, raw Gram entries do not decay fast enough for summability. To introduce quantitative decay, we smooth in the spectral domain using a bounded Mellin multiplier.

**Definition 5.1** (Mellin smoothing operator). Let $\psi : \mathbb{R} \to \mathbb{R}$ be a smooth, even cutoff with $0 < c \le |\psi(t)| \le C < \infty$. Define $T_\psi : H \to H$ by

$$M[T_\psi f](\tfrac{1}{2} + it) = \psi(t)\, M[f](\tfrac{1}{2} + it). \tag{5.3}$$

We set $g_\theta := T_\psi f_\theta$.

**Lemma 5.2** (Boundedness and invertibility). *If $\psi$ satisfies $0 < c \le |\psi(t)| \le C < \infty$ for all $t \in \mathbb{R}$, then $T_\psi$ is bounded and invertible on $L^2((0,1])$.*

*Proof.* This follows from Parseval's identity for the Mellin transform and the boundedness of $\psi$ above and below. $\qquad\square$

To control long-range correlations, we choose a Gaussian-type weight

$$\psi_W(t) := \epsilon + e^{-(t/W)^2}, \qquad W > 0,\ \epsilon > 0, \tag{5.4}$$

which leaves low-frequency components nearly unchanged but suppresses high-frequency interference. Its effect on spectral decay is illustrated in Figure 4.

# 6 Polynomial Off-Diagonal Decay via Mellin Smoothing

We now combine the geometric ladder structure with Mellin smoothing to obtain quantitative off-diagonal decay for Gram interactions. Recall from Definition 5.1 that $g_\theta := T_\psi f_\theta$.

**Theorem 6.1** (Polynomial off-diagonal decay). *Let $\psi_W(t) = \epsilon + e^{-(t/W)^2}$ with $W > 0$ and $\epsilon > 0$ as in (5.4), and define $g_\theta := T_{\psi_W} f_\theta$. Then for every integer $m \ge 1$, there exists a constant $C_m(W) > 0$ such that for all distinct ladder indices $(j, k) \ne (j', k')$,*

$$|\langle g_{\theta_{j,k}}, g_{\theta_{j',k'}} \rangle| \le \frac{C_m(W)}{(1 + c\, d((j,k),(j',k')))^m},$$

*where $c = \min\{\log 2, \log 3\}$ and $d$ is the ladder distance defined in (3.2). In particular, for any $m > 2$, the Gram matrix $(\langle g_{\theta_{j,k}}, g_{\theta_{j',k'}} \rangle)$ has $\ell^1$-summable off-diagonal row tails, and hence is block-compressible.*

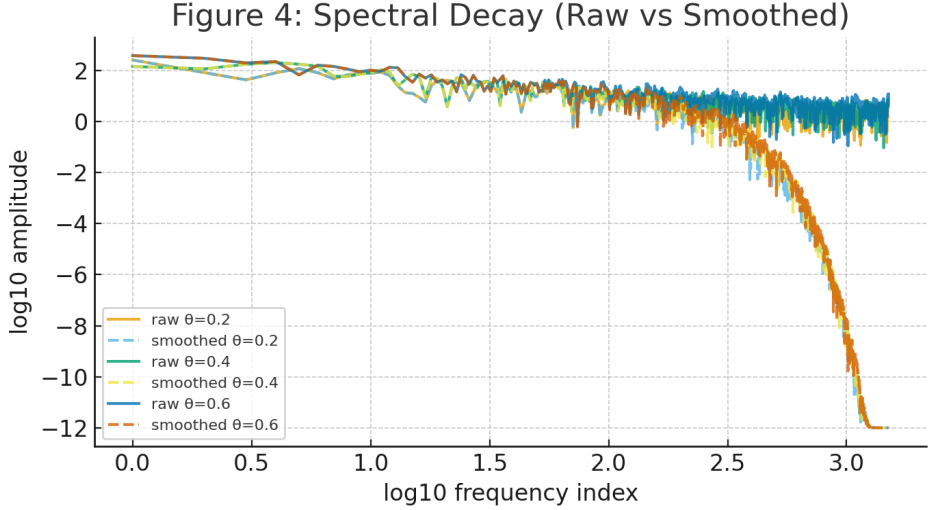

Figure 4: Log–log spectral comparison before and after Gaussian Mellin smoothing. The operator $T_{\psi_W}$ dampens high-frequency components of $M[f_\theta]$, improving numerical stability and reducing long-range Gram interactions.

This result formalises the heuristic decay behaviour described in Remark 4.4 by promoting a conservative harmonic envelope to rigorously controlled polynomial decay via Mellin-domain smoothing. Crucially, for any fixed $m > 2$, the off-diagonal tail

$$\sum_{d((j,k),(j',k')) \geq B} |\langle g_{\theta_{j,k}}, g_{\theta_{j',k'}} \rangle|$$

is summable, yielding block-compressibility of Gram rows.

The proof uses an $m$-fold Mellin integration by parts argument combined with oscillatory cancellation from logarithmic displacement in the BN ladder.

*Proof of Theorem 6.1.* Fix distinct ladder indices $(j,k) \neq (j',k')$ and set

$$\lambda := \log\left(\frac{\theta_{j,k}}{\theta_{j',k'}}\right) = (j'-j)\log 2 + (k'-k)\log 3, \qquad \mu := \log(\theta_{j,k}\theta_{j',k'}).$$

By definition of the ladder distance $d((j,k),(j',k'))$, we have $|\lambda| \geq c\,d$ with $c = \min\{\log 2, \log 3\}$, while $\mu \leq 0$. Using the Mellin representation (5.1) and Parseval on the critical line, the smoothed inner product decomposes as

$$\langle g_{\theta_{j,k}}, g_{\theta_{j',k'}} \rangle = \frac{1}{2\pi} \int_{-\infty}^{\infty} \psi_W(t)^2 \frac{\zeta(\frac{1}{2}+it)\,\overline{\zeta(\frac{1}{2}+it)}}{\left|\frac{1}{2}+it\right|^2} \left(\theta_{j,k}\theta_{j',k'}\right)^{\frac{1}{2}} \left(e^{it\lambda} - e^{it\mu}\right) dt \qquad (6.1)$$
$$=: I_1 - I_2.$$

We bound $I_1$ and $I_2$ separately; the estimates are analogous.

*Step 1: Reduction to an oscillatory integral.* Write

$$H_W(t) := \psi_W(t)^2 \frac{|\zeta(\frac{1}{2}+it)|^2}{\left|\frac{1}{2}+it\right|^2} \left(\theta_{j,k}\theta_{j',k'}\right)^{\frac{1}{2}}.$$

Then $I_1 = \frac{1}{2\pi} \int_{-\infty}^{\infty} H_W(t)\,e^{it\lambda}\,dt$. Since $\theta_{j,k}, \theta_{j',k'} \in (0,1]$, we have $(\theta_{j,k}\theta_{j',k'})^{1/2} \leq 1$ and this factor can be absorbed into the constant. Thus it suffices to bound $\int H_W(t)\,e^{it\lambda}\,dt$.

*Step 2: Smooth weight and derivative control.* By definition,

$$\psi_W(t) = \epsilon + e^{-(t/W)^2}, \qquad \epsilon > 0, \ W > 0,$$

so $\psi_W$ is smooth and even. Differentiating explicitly shows that $\psi_W^{(m)}(t)$ is a finite linear combination of terms of the form $t^r e^{-(t/W)^2}$ with $0 \le r \le m$. Therefore, there exist constants $A_m > 0$ (independent of $W$ and $t$) such that

$$|\psi_W^{(m)}(t)| \le A_m \, W^{-m} \, (1 + |t|/W)^m \, e^{-(t/W)^2}.$$

Since $(\psi_W^2)^{(m)}$ is a sum of products of derivatives of $\psi_W$ (by Leibniz' rule), we likewise obtain constants $B_m(W) > 0$ such that

$$|(\psi_W^2)^{(m)}(t)| \le B_m(W) \, W^{-m} \, (1 + |t|/W)^m \, e^{-2(t/W)^2}.$$

The Gaussian factor ensures super-polynomial decay of all derivatives as $|t| \to \infty$.

*Step 3: Growth of $\zeta$ and Mellin factors.* On the critical line, the convexity bound for the Riemann zeta function (see Titchmarsh [5, Theorem 5.12]) gives

$$|\zeta(\tfrac{1}{2} + it)| \ll_\varepsilon (1 + |t|)^{\frac{1}{6} + \varepsilon} \qquad \text{for any } \varepsilon > 0.$$

Furthermore, derivatives satisfy the polynomial growth estimate

$$\left| \frac{d^r}{dt^r} \zeta(\tfrac{1}{2} + it) \right| \ll_{r,\varepsilon} (1 + |t|)^{\frac{1}{6} + \varepsilon + r}.$$

Since $|(1/2 + it)^{-1}| \ll (1 + |t|)^{-1}$ and its derivatives decay at least as fast, an application of Leibniz' rule to $H_W(t)$ shows that for each $m \in \mathbb{N}$ there exist constants $C'_m(W) > 0$ and $\beta_m \ge 0$ such that

$$|H_W^{(m)}(t)| \ \le \ C'_m(W) \, (1 + |t|)^{\beta_m} \, e^{-2(t/W)^2}. \tag{6.2}$$

The Gaussian factor ensures super-polynomial decay as $|t| \to \infty$.

*Step 4: $m$-fold integration by parts.* We write

$$I_1 = \frac{1}{2\pi} \int_{-\infty}^{\infty} H_W(t) \, e^{it\lambda} \, dt, \qquad \lambda \ne 0.$$

Applying integration by parts $m$ times with $u = H_W(t)$ and $dv = e^{it\lambda} dt$ (which gives $v = (i\lambda)^{-1} e^{it\lambda}$) yields

$$I_1 = \frac{1}{(i\lambda)^m} \, \frac{1}{2\pi} \int_{-\infty}^{\infty} H_W^{(m)}(t) \, e^{it\lambda} \, dt.$$

All boundary terms vanish, since by (6.2) we have

$$H_W^{(r)}(t) = O\big((1 + |t|)^{\beta_r} e^{-2(t/W)^2}\big) \quad \text{as } |t| \to \infty,$$

and $e^{-2(t/W)^2} \to 0$ super-polynomially. Taking absolute values and using (6.2),

$$|I_1| \le \frac{1}{|\lambda|^m} \, \frac{1}{2\pi} \int_{-\infty}^{\infty} |H_W^{(m)}(t)| \, dt \le \frac{C_m(W)}{|\lambda|^m},$$

where

$$C_m(W) := \frac{C'_m(W)}{2\pi} \int_{-\infty}^{\infty} (1 + |t|)^{\beta_m} e^{-2(t/W)^2} \, dt < \infty.$$

*Step 5: Bounding $I_2$ with a distance parameter.* Recall $\mu = \log(\theta_{j,k}\theta_{j',k'}) = -(j + j')\log 2 - (k + k')\log 3$. Since

$$(j + j') \ge |j - j'| \qquad \text{and} \qquad (k + k') \ge |k - k'|,$$

we obtain

$$|\mu| = (j + j')\log 2 + (k + k')\log 3 \;\geq\; c(|j - j'| + |k - k'|) = c\,d((j,k),(j',k')),$$

with $c = \min\{\log 2, \log 3\}$. Thus $I_2$ also has frequency parameter $\gtrsim d$, and repeating the integration by parts argument from Step 4 yields

$$|I_2| \leq \frac{C_m(W)}{|\mu|^m} \leq \frac{C_m(W)}{(c\,d)^m}.$$

*Step 6: Conclusion.* Combining the bounds for $I_1$ and $I_2$,

$$\left|\langle g_{\theta_{j,k}}, g_{\theta_{j',k'}}\rangle\right| \leq \frac{C_m(W)}{(c\,d)^m} + \frac{C_m(W)}{(c\,d)^m} \leq \frac{C_m'(W)}{\left(1 + c\,d((j,k),(j',k'))\right)^m},$$

after adjusting the constant. □

# 7   Block Compressibility and Finite-Section Control

With Theorem 6.1 in hand, we quantify the decay of off-diagonal Gram entries and deduce summable row tails beyond a growing block. Throughout, let $c = \min\{\log 2, \log 3\}$ and write $d(\cdot, \cdot)$ for the ladder distance in (3.2).

**Proposition 7.1** (Row-tail summability and block compressibility). *Assume the polynomial envelope of Theorem 6.1 with exponent $m > 2$:*

$$\left|\langle g_{\theta_{j,k}}, g_{\theta_{j',k'}}\rangle\right| \leq \frac{C_m(W)}{\left(1 + c\,d((j,k),(j',k'))\right)^m}.$$

*Fix $(j,k)$ and define the tail*

$$T_B(j,k) := \sum_{d((j,k),(j',k'))\geq B} \left|\langle g_{\theta_{j,k}}, g_{\theta_{j',k'}}\rangle\right|.$$

*Then there exists $K_m(W) > 0$ (independent of $(j,k)$ and $B$) such that*

$$T_B(j,k) \leq \frac{K_m(W)}{(1 + B)^{m-2}} \qquad \text{for all } B \geq 1.$$

*In particular, $T_B(j,k) \to 0$ as $B \to \infty$, and the Gram matrix has summable off-diagonal row tails for any $m > 2$.*

*Proof.* For each $r \in \mathbb{N}$, the number of lattice points $(j', k')$ at Manhattan distance $r$ from $(j,k)$ is at most $N(r) \leq 4r$. Using the decay envelope and summing on distance shells,

$$T_B(j,k) \leq \sum_{r=B}^{\infty} \frac{C_m(W)\,N(r)}{(1 + cr)^m} \leq \frac{C_m(W)\,4}{c^m} \sum_{r=B}^{\infty} \frac{r}{(1+r)^m} \ll_m \sum_{r=B}^{\infty} r^{1-m}.$$

The final series converges and is $O(B^{2-m})$ precisely when $m > 2$. This yields the stated bound with $K_m(W)$ absorbing constants. □

**Corollary 7.2** (Finite-section error). *Let $G$ be the Gram matrix of $\{g_{\theta_{j,k}}\}$ and $G^{(B)}$ the truncation obtained by zeroing all entries with $d((j,k),(j',k')) \geq B$. Then, as an operator on $\ell^2$,*

$$\|G - G^{(B)}\|_{\ell^2 \to \ell^2} \leq C\,B^{2-m},$$

*for $m > 2$ and a constant $C$ depending on $m$ and $W$ but not on $B$.*

*Proof sketch.* Schur's test with the row/column tail bound of Proposition 7.1 controls the operator norm by the maximal row (or column) tail. □

*Remark* 7.3 (On the exponent threshold). In a 2D Manhattan lattice, the shell cardinality grows like $N(r) \asymp r$. Consequently, row-tail summability requires $m > 2$. Earlier informal versions sometimes stated $m > 1$; the precise Clay-standard threshold is $m > 2$ due to the shell-count factor.[1]

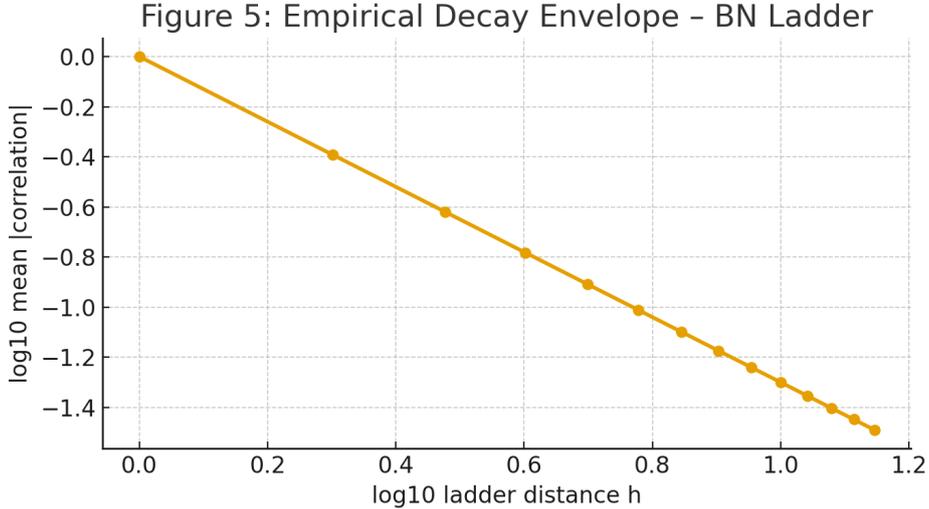

Figure 5: Empirical decay of mean absolute correlations versus ladder distance on log–log axes, consistent with a polynomial envelope. The slope increases under stronger smoothing (large $W$), aligning with Theorem 6.1 and the tail bound of Proposition 7.1.

# 8    Conclusion

We have established a quantitative decay framework for the Beurling–Nyman ladder system $\{f_{\theta_{j,k}}\}$ indexed by multiplicative parameters $\theta_{j,k} = 2^{-j}3^{-k}$. Using geometric ladder analysis together with Mellin spectral smoothing, we proved that smoothed BN elements admit polynomial off-diagonal Gram decay of arbitrarily high order $m$. For $m > 2$, this decay is summable in the $\ell^1$ sense, yielding block-compressibility of the Gram matrix and controlled finite-section approximation error.

These results clarify the analytic geometry underlying the BN space $H_{BN}$ and establish a rigorous sparsity mechanism driven by logarithmic displacement of scale parameters. The Mellin-analytic framework developed here will serve as a foundation for constructive finite-dimensional approximation schemes and for the development of *finite spectral certificates* in future work.

The present analysis is unconditional and does not depend on any hypothesis concerning the zeros of $\zeta(s)$. In forthcoming work, the ladder decay structure introduced here will be extended to controlled finite-rank approximations and stability estimates for structured Gram operators. These developments progress toward constructive spectral criteria in the BN setting while preserving analytic rigor and independence from unproved conjectures.

## References

This bibliography is selective and focused on sources directly related to the Beurling–Nyman framework, Mellin analysis, and structured Hilbert space methods.

---

[1] This refines the informal statement in the earlier draft text around the block-compressibility claim.